\documentclass{article}

\title{Seshadri constants in finite subgroups of abelian surfaces.}
\author{Luis Fuentes Garc\'{\i}a}
\date{}
\newtheorem{teo}{Theorem}[section]

\newtheorem{cor}[teo]{Corollary}
\newtheorem{lemma}[teo]{Lemma}

\newtheorem{rem}[teo]{Remark}
\newtheorem{conj}[teo]{Conjecture}

\def\g2{\pi}

\newcommand\Te{{\cal O}}

\def\ZZ{\leavevmode\hbox{$\rm Z$}}

\def\qed{\hspace{\fill}$\rule{2mm}{2mm}$}

\newcommand\lrw{\longrightarrow}

\begin{document}

\maketitle

\begin{center}
{\it \footnotesize\small Luis Garc\'{\i}a L\'opez in memoriam. \normalsize}
\end{center}

\vspace{0.1cm}

\begin{abstract}

Given an \'etale quotient $q:X\lrw Y$ of smooth projective varieties we relate the simple Seshadri constant of a line bundle $M$ on $Y$ with the multiple Seshadri constant of $q^*M$ in the points of the fiber. We apply this method to compute the Seshadri constant of polarized abelian surfaces in the points of a finite subgroup.



\end{abstract}

\section{Introduction.}

The multiple Seshadri constants are a natural generalization of the Seshadri constants at single points defined by Demailly in \cite{De92}. If $X$ is a smooth projective variety of dimension $n$, $L$ is an ample line bundle on $X$ and $x_1,\ldots,x_r$ are distinct points in $X$, then the Seshadri constant of $L$ at $x_1,\ldots,x_r$ is:
$$
\epsilon(L;x_1,\ldots,x_r)=sup\{\epsilon|\, f^*L-\epsilon \sum_{i=1}^rE_i\hbox{ is nef }\},
$$
where $f$ is the blowing up of $X$ at $x_1,\ldots,x_r$  and $E_1,\ldots,E_r$ are the exceptional divisors. These constants have the upper bound:
$$
\epsilon(L; r)\leq \sqrt[n]{\frac{L^{n}}{r}}.
$$
However, explicit values are difficult to obtain even when $r=1$. General bounds for the simple Seshadri constants on surfaces are given in \cite{Ba99}, \cite{Na03} or \cite{St98}. They were computed for simple abelian surfaces by Th. Bauer (see \cite{Ba99}); Ch. Schultz gave values for Seshadri constants on products of two elliptic curves (see \cite{Sc04}).

The case of the multiple Seshadri constants is harder. For example, in the plane the Nagata conjecture is still an open problem (see \cite{StSz04}):

\begin{conj}[Nagata conjecture]

Let $x_1,\ldots,x_r$ be $r\geq 10$ be general points in $P^2$  then:
$$
\epsilon(\Te_{P^2}(1); x_1,\ldots,x_r)=\frac{1}{\sqrt{r}}.
$$
\end{conj}

This has been extended for an arbitrary surface. When $r$ is big enough, the value of the Seshadri constant at $r$ very general points is conjectured to be maximal (see \cite{Bi99}). Very intersting lower bounds for multiple Seshadri constants were given by B. Harbourne in \cite{Ha03}. In \cite{Tu-Ga04} Tutaj-Gasi\'nska give bounds for the Seshadri constant of abelian surfaces in half-periods points; in \cite{Tu-Ga05}, he gives the exact values in two half-periods points.

In this paper, we obtain the exact value of the multiple Seshadri constants of polarized abelian surfaces in points of a finite subgroup. This generalizes the results of \cite{Tu-Ga04} and \cite{Tu-Ga05} but applying a different method.

If $q:X\lrw Y$ is an \'etale quotient of smooth projective varieties we prove that the simple Seshadri constant of a line bundle $M$ on $Y$ is the same that the multiple Seshadri constant of $q^*M$ in the points of the fiber. We apply this result when $X$ is an abelian surface and $Y=X/G$ is the quotient by a finite subgroup $G$. Since, the  simple Seshadri constants on abelian surfaces are known (see \cite{Ba99}, \cite{Sc04}), we obtain the multiple Seshadri constants on $G$.

In particular, when $X$ is an abelian surface with Picard number one, we prove the following:

\begin{teo}\label{principala}

Let $(X,L)$ be a polarized abelian surface of type $(1,d)$ with $\rho(X)=1$. Let $x$ be a point of $X$. Let $G$ be a finite subgroup of $X$ of order $g$. Consider the \'etale quotient:
$$
q:X\lrw X/G
$$
Let $n$ be the minor integer verifying $nL=q^*M$ for some line bundle $M$ on $X/G$. Then:

\begin{enumerate}

\item If $\sqrt{2d/g}$ is rational, then $\epsilon(L; x+G)=\sqrt{\frac{2d}{g}}$.

\item If $\sqrt{2d/g}$ is irrational, then 
$$
\epsilon(L; x+G)= \frac{k_0}{l_0}\frac{2dn}{g}=\sqrt{1-\frac{1}{l_0^2}}\sqrt{\frac{L^2}{g}}
$$
where $(l_0,k_0)$ is the primitive solution of Pell's equation $l^2-\frac{2n^2d}{g}k^2=1$.

\end{enumerate}

\end{teo}

Note, that the order of $G$ and the degree of $L$ is not sufficient to determine the value of the Seshadri constants. It depends also on the the structure of $G$ (see Remark \ref{menor} for details).

In general, we also prove:

\begin{teo}

Let $(X,L)$ be a polarized abelian surface. The multiple Seshadri constant of $L$ at the points of a finite subgroup is rational.

\end{teo}

\section{The main theorem.}

\begin{teo}\label{principal}

Let $q:X\lrw Y$ be an \'etale $n:1$ quotient between to smooth projective varieties. Let $M$ be a line bundle on $Y$ and $y\in Y$. Then:
$$
\epsilon(q^*M; q^{-1}(y))=\epsilon(M,y)
$$

\end{teo}
{\bf Proof:} Let $g:\tilde{Y}\lrw Y$ be the blowing up of $Y$ at $y$. Let $f:\tilde{X}\lrw X$ be the blowing up of $X$ at $q^{-1}(y)$. There is an induced morphism $\tilde{q}:\tilde{X}\lrw \tilde{Y}$ such that the following diagram is commutative:
$$
\matrix{
\phantom{f}\tilde{X} & \stackrel{\tilde{q}}{\lrw} & \tilde{Y}\phantom{g} \cr
 f\downarrow & & \downarrow g \cr
\phantom{f}X & \stackrel{q}{\lrw} & Y\phantom{g} }
$$
If $E$ is the exceptional divisor of $g$ and $E_1,\ldots,E_n$ are the exceptional divisors of $f$, we have:
$$
f^*q^*M-\epsilon(E_1+\ldots+E_n)=\tilde{q}^*g^*M-\epsilon\tilde{q}^*E=\tilde{q}^*(g^*M-\epsilon E).
$$
Thus,
$$
f^*q^*M-\epsilon(E_1+\ldots+E_n)\hbox{ is nef }\iff g^*M-\epsilon E\hbox{ is nef}.
$$ \qed

\section{Multiple Seshadri constant on polarized abe\-lian surfaces.}

We will apply the main theorem to compute the multiple Seshadri constants on polarized abelian surfaces. We will use the results of Bauer to compute the simple Seshadri constant (see \cite{Ba99}):

\begin{teo}[Bauer]\label{bauer}

Let $(Y,M)$ be an abelian surface of type $(1,d')$ with $\rho(Y)=1$, 

\begin{enumerate}

\item If $\sqrt{2d'}$ is rational, then $\epsilon(M)=\sqrt{2d'}.$

\item If $\sqrt{2d'}$ is irrational, then
$$
\epsilon(M)= \frac{k_0}{l_0}2d'=\sqrt{1-\frac{1}{l_0^2}}\sqrt{M^2}
$$
where $(k_0,l_0)$ is the primitive equation of the Pell's equation $l^2-2dk^2=1$.

\end{enumerate}

\end{teo}

Let us recall some basic facts about abelian surfaces. We will follow the notation of \cite{BiLa92}. Let $(X,L)$ be a polarized abelian surface of type $(1,d)$. Let $H$ be the first Chern class of $L$. There is a basis $\lambda_1,\lambda_2,\mu_1,\mu_2$, respect to which $E=Im H$ is given by the matrix:
$$
\left(\matrix{\phantom{-}0 & D \cr -D & 0}\right)
$$
where $D=diag(1,d)$. In this way, the abelian surface $X$ is the quotient 
$$
\pi:V\lrw X=V/\Lambda
$$
with
$$
\Lambda=\langle \lambda_1,\lambda_2\rangle \oplus \langle \mu_1,\mu_2 \rangle.
$$

Let $G$ be a finite subgroup of $X$ of order $g$. Consider the quotient map:
$$
q:X\lrw X/G.
$$
The variety $Y=X/G$ is an abelian surface. In particular $Y=V/\Lambda'$, where $\Lambda'=\pi^{-1}(G)$. There is a criterion for a line bundle $L'\in Pic(X)$ to descend under $q$:

\begin{lemma}

Let $L'=L(H',\chi')$ be a line bundle on $X$. Then, $L'=q^*M$ for some line bundle $M\in Pic(Y)$ if and only if $Im H'(\Lambda',\Lambda')\subset \ZZ$. 

\end{lemma}
{\bf Proof:} See Chapter 2, Corollary 4.4 of \cite{BiLa92}. \qed

\begin{cor}

If $L$ is a line bundle on $X$ of type $(1,d)$, then $exp(G)^2L=q^*M$ for some $M\in Pic(Y)$.

\end{cor}
{\bf Proof:} It is a consequence of the previous corollary. The exponent of a group $G$ is the least common multiple of the orders of the elements of $G$. The first Chern class of $exp(G)^2L$ is $exp(G)^2H$. Let $x,x'\in \pi^{-1}(G)$. We know that $\exp(G)\cdot x,\exp(G)\cdot x'\in \Lambda$. Thus:
$$
Im\, (exp(G)^2H)(x,x')=Im\, H(exp(G)\cdot x,exp(G)\cdot x')\in \ZZ.
$$ \qed

\begin{rem}\label{menor}

We can consider the minor integer $n$ such that $nL$ descend to a line bundle $M$. The minimality of $n$ implies that $M$ is a primitive line bundle on $Y$. We know that $1\leq n\leq exp(G)^2$. However, the number $n$ does not depend only on the exponent of $G$. For example, if $G$ is the cyclic group generated by $\pi(\lambda_2/d)$ then $L$ descends to a line bundle $M$ of type $(1,1)$, so $n=1$. On the other hand, if $G$ is the  group generated by $\pi(\lambda_1/k),\pi(\mu_1/k)$ for any $k>1$ then the minor value of $n$ is $n=k^2=exp(G)^2$.

\end{rem}

The main result on abelian surfaces with Picard number one will be the following:

\begin{teo}

Let $(X,L)$ be a polarized abelian surface of type $(1,d)$ with $\rho(X)=1$. Let $x$ be a point of $X$. Let $G$ be a finite subgroup of $X$ of order $g$. Consider the \'etale quotient:
$$
q:X\lrw X/G
$$
Let $n$ be the minor integer verifying $nL=q^*M$ for some line bundle $M$ on $X/G$. Then:

\begin{enumerate}

\item If $\sqrt{2d/g}$ is rational, then $\epsilon(L; x+G)=\sqrt{\frac{2d}{g}}$.

\item If $\sqrt{2d/g}$ is irrational, then 
$$
\epsilon(L; x+G)= \frac{k_0}{l_0}\frac{2dn}{g}=\sqrt{1-\frac{1}{l_0^2}}\sqrt{\frac{L^2}{g}}
$$
where $(l_0,k_0)$ is the primitive solution of Pell's equation $l^2-\frac{2n^2d}{g}k^2=1$.

\end{enumerate}

\end{teo}
{\bf Proof:} Let $y=q(x)$. By the Theorem \ref{principal}:
$$
\epsilon(L; x+G)=\frac{1}{n}\epsilon(nL; x+G)=\frac{1}{n}\epsilon(M,y).
$$
The line bundle $M$ is a primitive line bundle of type $(1,d')$ with $d'=n^2d/g$. Thus $\sqrt{2d/g}$ is rational if and only if $\sqrt{2d'}$ is rational. Now, the result follows from Theorem \ref{bauer}. \qed

\begin{cor}

Let $(X,L)$ be a polarized abelian surface of type $(1,d)$ with $\rho(X)=1$. Let $x_1,\ldots,x_r$ be $r$  general points of $X$. Then:

\begin{enumerate}

\item If $\sqrt{2d/r}$ is rational, then $\epsilon(L; x_1,\ldots,x_r)=\sqrt{\frac{2d}{r}}$.

\item If $\sqrt{2d/r}$ is irrational, then 
$$
\epsilon(L; x_1,\ldots,x_r)\geq 2d\frac{k_0}{l_0}=\sqrt{1-\frac{1}{l_0^2}}\sqrt{\frac{L^2}{r}}
$$
where $(l_0,k_0)$ is the primitive solution of Pell's equation $l^2-2rdk^2=1$.

\end{enumerate}

\end{cor}
{\bf Proof:} By the semicontinuity of the Seshadri constant,
$$
\epsilon(L; x_1,\ldots,x_r)\geq \epsilon(L; x+G)
$$
for any point $x\in X$ and any subgroup $G$ of order $r$. In particular, taking the cyclic subgroup $G=\langle \pi(\lambda_1/r) \rangle$ and applying the previous theorem we obtain the desired bound. \qed

\begin{cor}

Let $(X,L)$ be a polarized abelian surface of type $(1,d)$ with $\rho(X)=1$. Let $x$ be a point of $X$.  Let $X_m$ be the subgroup of $m$-torsion points. Suppose that $\sqrt{2d^2}$ is not an integer, then:
$$
\epsilon(L; x+X_m)= 2\frac{d}{m^2}\frac{k_0}{l_0}=\sqrt{1-\frac{1}{l_0^2}}\sqrt{\frac{L^2}{m^4}}
$$
where $(l_0,k_0)$ is the primitive solution of Pell's equation $l^2-2dk^2=1$.

\end{cor}
{\bf Proof:} Note that, in this case, the minor number $n$ such that $nL$ descends under $q$ is $n=m^2$. Now, it is sufficient to apply the Theorem \ref{principala}. \qed

\begin{rem}

In \cite{Tu-Ga04}, Tutaj-Gasi\'nska obtains a bound for the Seshadri constant in half-periods of a line bundle on a polarized  abelian surface of type $(1,d)$:
$$
\epsilon(L; X_2)\leq 2\sqrt{1-\frac{1}{l_0^2}}\sqrt{\frac{L^2}{16}}
$$
where $(l_0,k_0)$ is the primitive solution of Pell's equation $l^2-32dk^2=1$. Here, we see that we exact value appears when we use the Pell's equation $l^2-2dk^2=1$. 

\end{rem}

\begin{cor}\label{doble}

Let $(X,L)$ be a polarized abelian surface of type $(1,d)$ with $\rho(X)=1$. Let $e_1,e_2$ be two half periods of $X$. Suppose that $\sqrt{d}$ is not an integer.

\begin{enumerate}

\item If $e_1-e_2\in K(L)$ then 
$$
\epsilon(L; e_1,e_2)=d\frac{k_0}{l_0}
$$
where $(l_0,k_0)$ is the primitive solution of Pell's equation $l^2-dk_0^2=1$.

\item If $e_1-e_2\not\in K(L)$ then
$$
\epsilon(L; e_1,e_2)=2d\frac{k_0}{l_0}
$$
where $(l_0,k_0)$ is the primitive solution of Pell's equation $l^2-4dk_0^2=1$.

\end{enumerate}

\end{cor}
{\bf Proof:} Consider the group $G=\langle e_1-e_2 \rangle$. It has order two. Moreover, $L$ descends to a line bundle in $X/G$ if and only if $e_1-e_2\in K(L)$. Now we only have to apply the Theorem \ref{principala}. \qed

\begin{rem}

In \cite{Tu-Ga05}, a similar study is made. However, the Theorem (1, \cite{Tu-Ga05}) does not distinguish which are the two half-periods where the Seshadri constant is computed. The problem is the application of Lemma (19, \cite{Tu-Ga05}). It says that one can choose a line bundle $L_c$ such that any two half periods  $e_1,e_2$ have the same parity. However, this changes the number of even and odd half periods of $L$. From this, the arguments on page (532, \cite{Tu-Ga05}) can fail.

\end{rem}

We could use similar arguments to compute the multiple Seshadri constants at points of a finite subgroup of line bundles on abelian surfaces with Picard number great than $1$. The simple Seshadri constants of these surfaces were computed by C. Schultz in \cite{Sc04}. Anyway, we know that the simple Seshadri constant of a line bundle on any abelian surface is always rational. From this:

\begin{teo}

Let $(X,L)$ be a polarized abelian surface. The multiple Seshadri constant of $L$ at the points of a finite subgroup is rational.

\end{teo}

{\bf E-mail:} lfuentes@udc.es

Luis Fuentes Garc\'{\i}a. 

Departamento de M\'etodos Matem\'aticos y Representaci\'on.

E.T.S. de Ingenieros de Caminos, Canales y Puertos. 

Universidad de A Coru\~na. Campus de Elvi\~na. 15192 A Coru\~na (SPAIN)


\begin{thebibliography}{77}

\bibitem{Ba98}{\sc Bauer, T. }
{\it Seshadri constants and periods of polarized abelian varieties. }
Math. Ann. {\bf 312}, 607-623 (1998).

\bibitem{Ba99}{\sc Bauer, T. }
{\it Seshadri constants on algebraic surfaces. }
Math. Ann. {\bf 313}, 547-583 (1999).

\bibitem{BiLa92}{\sc Birkenhake, Ch.; Lange, H.}
{\it Complex abelian vareties.} Springer-Verlag (1992).

\bibitem{Bi99}{\sc Biran, P.}
{\it Constructing new ample divisors out of old ones.}
Duke Math. J. {\bf 98} 113-135, (1999).

\bibitem{De92}{\sc Demailly, J.-P.}
{\it Singular hermitian metrics on positive line bundles}
Lecture Notes Math. {\bf 1507}, 87-104 (1992).

\bibitem{Ha03}{\sc Harbourne, B.}
{\it Seshadri constants and very ample divisors on algebraic surfaces.}
J. Reine Angew. Math {\bf 559}, 115-122 (2003).	

\bibitem{Na03}{\sc Nakamaye, M.J.}
{\it Seshadri constants and the geometry of surfaces.}
J. Reine Angew. Math. {\bf 564}, 205-214 (2003).

\bibitem{Sc04}{\sc Schultz, Ch.}
{\it  Seshadri constants on abelian surfaces.}
Thesis, Marburg, 2004.

\bibitem{St98}{\sc Steffens, A.} 
{\it Remarks on Seshadri constants.}
Math. Z. {\bf 227}, 505-510 (1998).

\bibitem{StSz04}{\sc Strycharz-Szemberg, B.; Szemberg, T.}
{\it  Remarks on the Nagata conjecture.}
Serdica Math. J. {\bf 30}, Nº 2-3, 405-430 (2004)

\bibitem{Tu-Ga04}{\sc	Tutaj-Gasi\'nska, H.}
{\it Seshadri constants in half-periods of an abelian surface.}
J. Pure Appl. Algebra {\bf 194}, Nº1-2, 183-191 (2004).

\bibitem{Tu-Ga05}{\sc	Tutaj-Gasi\'nska, H.}
{\it Seshadri constants in two half periods.}
Arch. Math. {\bf 85}, Nº. 6, 514-526 (2005).







\end{thebibliography}
\end{document}